\documentclass[graybox]{svmult}

\usepackage{mathptmx}       
\usepackage{helvet}         
\usepackage{courier}        
\usepackage{type1cm}        
\usepackage{makeidx}         
\usepackage{graphicx}        
\usepackage{multicol}        
\usepackage[bottom]{footmisc}

\usepackage{amsfonts,amsmath}

\makeindex             

\begin{document}

\title*{Folding patterns in partially delaminated thin films}
\author{David Bourne, Sergio Conti and Stefan M\"uller}
\institute{
David Bourne \at
 Department of Mathematical Sciences, Durham University
 \and
 Sergio Conti \at
 Institut f\"ur Angewandte Mathematik,
Universit\"at Bonn
\and
Stefan M\"uller \at
 Institut f\"ur Angewandte Mathematik,
Universit\"at Bonn}
%
%
\maketitle

\abstract{Michael Ortiz and Gustavo Gioia showed in the 90s
that the complex patterns arising in compressed
elastic films
can be analyzed within the context of the
calculus of variations. Their initial work focused
on films partially debonded from
the substrate, subject to isotropic compression arising
from the difference in thermal expansion coefficients between film and substrate.
In the following two decades different
geometries have been studied, as for example anisotropic compression.
We review recent mathematical progress in this area,
focusing on the rich phase diagram of partially debonded films with a
lateral boundary condition.}

\section{Introduction}

Elastic films deposited on a substrate are often subject, after
thermal expansion, to compressive strains which are released
by debonding and buckling, generating a variety of microstructures.
The work of Michael Ortiz and Gustavo Gioia in the 90s {[1, 2]} opened
the way for the use of the tools of
calculus of variations in the study of these structures.
Their starting point was the  F\"oppl-von K\'arm\'an plate theory, as given in ({4})
below.
One of their insights was that the key
nonconvexity which gives rise to the microstructure can be understood in terms
of the out-of-plane displacement alone, leading after some rescalings
to the Eikonal functional, as given in
({1}) below.
This functional contains a term of the form $(|D w|^2-1)^2$, where $w$ is the normal displacement, which favours
deformations with the property that the gradient of $w$ is approximately a unit vector, independently of the orientation. Since the
film is still bound to the substrate at the boundary of the debonded region, the appropriate boundary
condition is $w=0$, which prescribes that the average over $\Omega$ of the gradient of $w$ vanishes. Therefore the resulting low-energy
deformations have gradient $D w$ oscillating between different values. As in many nonconvex variational problems,
oscillations on very small scales may be energetically convenient,
see {[3, 4]}.
Correspondingly, the variational problem $\int_\Omega (|Dw|^2-1)^2dx$ is not
lower semicontinuous, and  - depending on the boundary data and forcing - does not have minimizers.
However, the curvature term  $\sigma^2 |D^ 2w|^2$ penalizes oscillations on an exceedingly fine scale and thereby ensures existence
of minimizers. The solutions then have oscillations on an intermediate scale, which is determined by the competition
between the two terms.
The analysis of the specific functional proposed by Ortiz and Gioia is reviewed
in Section {2} below.

The approach of Ortiz an Gioia was later extended to the full vectorial  F\"oppl-von K\'arm\'an
energy, and also to three-dimensional elasticity. These refinements 
explained the
appearance of oscillations on two different length scales, with
coarse oscillations in a direction normal to the boundary, and fine oscillations in the
direction tangential to the boundary, as discussed in Section {3} below.

Recently interest has been directed to controlling the microstructures by
designing the geometry of the debonded region appropriately {[5, 6]}.  The key idea is to introduce a sacrificial layer between
the film and the substrate, and then to selectively etch away a part of it, so that the boundary of the debonded region is straight. The film then partially rebonds to the surface, leading to complex patterns of tunnels.
A study of these patterns within the Ortiz-Gioia framework,
with a variational functional containing the F\"oppl-von K\'arm\'an   energy and a fracture term
 proportional to the debonded area, is presented in Section {4}.
The mathematical analysis leading to the upper bounds of Theorem {6} suggests the presence of different
types of patterns in different parameter ranges.
The picture is rather easy in the two extreme cases in which the bonding energy per unit area is
very small or very large.
Indeed, in the first one the patterns observed for completely debonded films
give the optimal energy scaling, in the second one the optimal state corresponds to the film
completely bound to the substrate. In the intermediate regime we expect a richer picture,
with bonded areas separated by thin debonded tunnels. For a certain regime, depending on the relation between the bonding
energy per unit area, the film thickness and the compression ratio, a construction in which the tunnels
branch and refine close to the boundary has a lower energy than the one with straight tunnels,
see discussion in  Section {4} below.
The microstructure formation in thin films can be understood at a qualitative level as a form of
Euler buckling instability. The relevant experiments, however, are well beyond the stability
threshold, as discussed in Section {5} below.

\begin{figure}[t]\centering
\includegraphics[width=0.55\linewidth]{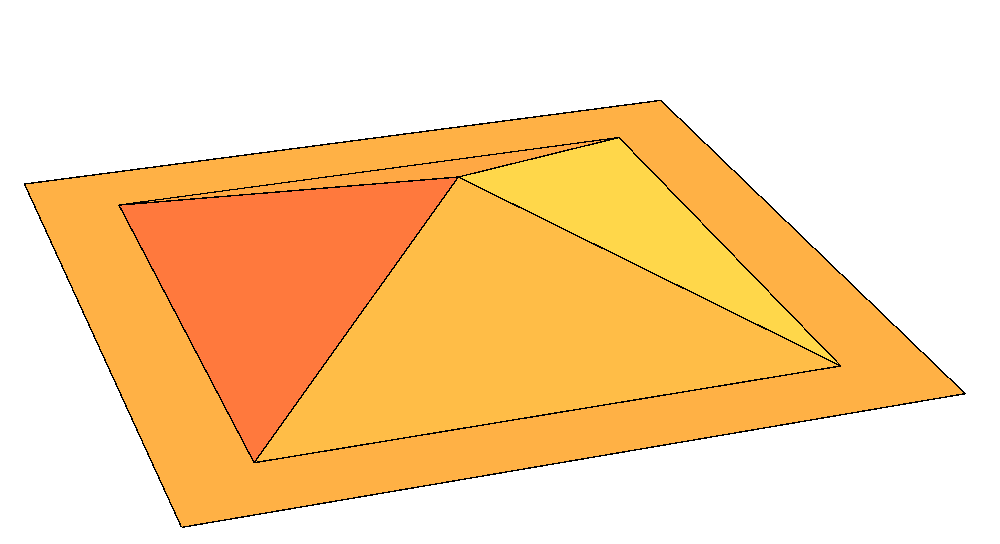}
\caption{Sketch of a deformation achieving the optimal energy in ({1}).
Here the debonded region $\Omega=(0,1)^2$ is a square, and
the distance to the boundary gives a ``tent''-form. The convolution in ({3})
makes the folds then smooth transitions on a small scale.}
\end{figure}

\section{Scalar modeling of compressed thin films}
Ortiz and Gioia showed that, if tangential displacements are neglected, the energy
of a compressed thin film can be characterized by the functional
\begin{equation}
  I_\sigma[w]=\int_\Omega  \left(\left(|Dw|^2-1\right)^2 + \sigma^2 |D^2w|^2 \right)dx,
\end{equation}
subject to $w=0$ on $\partial\Omega$
and $w\ge 0$ in $\Omega$.
Here $\Omega\subset\mathbb{R}^2$
represents the debonded region, $w:\Omega\to[0,\infty)$ the rescaled normal displacement, and  $\sigma$ is a small
parameter related to the thickness of the film.
This functional arises also naturally in the study of liquid crystal configurations
{[7]} and of magnetic structures in thin films {[8]}.
Despite a large mathematical effort
{[7, 9, 10, 11, 12, 13, 14, 15]} the
problem ({1}) is not yet completely
understood; it has been shown that the minimal energy
is proportional to $\sigma$ but the $\Gamma$-limit of $\sigma^{-1}I_\sigma[w]$ has only been partially identified.
The natural candidate is
\begin{equation}
 I_0[w]=\frac13\int_{J_{D w}} | [D w]|^3 d\mathcal{H}^1
\end{equation}
restricted to functions $w:\Omega\to \mathbb{R}$ which solve the Eikonal equation $|D w|=1$ and are sufficiently regular.
Here, $J_{Dw}$ denotes the set of points (typically, a curve) where the gradient $Dw$ is not continuous,
$[Dw]$ denotes its jump  across the interface, and $d\mathcal{H}^1$ the line integral along the interface.
In particular, under the additional assumption that  $Dw$ is a function of bounded variation, it has
been shown that for $\sigma\to0$ the scaled functionals $\sigma^{-1}I_\sigma$ converge,
in the sense of $\Gamma$-convergence, to $I_0$, see  {[11, 12, 13]} for the lower bound and {[14, 15, 16]} for the upper bound. However, it is also clear that finiteness
of the energy does not imply that $Dw$ has bounded variation, but only that $w$ belongs to
a larger space, called $AG(\Omega)$, see  {[11, 10]}. Therefore the result is still incomplete.

The Eikonal equation $|Dw|=1$ with the boundary data $w=0$ on $\partial\Omega$ is solved by the distance to the
boundary, $w_0(x)=\mathrm{dist}(x,\partial\Omega)$. For example, if $\Omega$ is a square this leads to the tent-shaped
deformation illustrated in Figure {1}. The function $w_0$ is however only Lipschitz continuous, not twice
differentiable, and makes the curvature
term $\int_\Omega \sigma^2 |D^2w|^2dx$ infinite. Therefore Ortiz and Gioia  {[1, 2]}  proposed to use a smoothed version
of the distance function,
\begin{equation}
 w_\sigma(x) = \int_\Omega \mathrm{dist}(y,\partial\Omega) \varphi_\sigma(x-y)dy
\end{equation}
where $\varphi_\sigma$ is a mollifier on the scale $\sigma$, i.e., $\varphi_\sigma\in C^\infty_c(B_\sigma)$
with $\int_{\mathbb{R}^2} \varphi_\sigma\,dx=1$ and $|D\varphi_\sigma|\le c/\sigma^3$. Then the regularized gradient $Dw_\sigma$ has length
close to 1 on most of the domain $\Omega$, but at boundaries between regions where $Dw_0$ has different orientations
$Dw_\sigma$ changes smoothly over a length scale $\sigma$ from one value to the other. The bending energy is correspondingly localized
in a stripe of thickness $2\sigma$ around the interfaces, see Figure {2}.
The prediction that
minimizers of ({1}) are well represented by  $w_\sigma$ is in good agreement,
at least for some geometries, with experimental observations {[1, 2]}.

The work of Ortiz and Gioia was then extended to related problems, showing for example that
under  anisotropic compression
 branching-type microstructures appear close to the boundary
{[17, 18]}, or that in certain regimes  telephone-cord blisters
develop {[19, 20, 21]} thanks to the interaction between the
elastic deformation and the fracture problem that determines the boundary of the debonded region.

\begin{figure}[t]\centering
\includegraphics[width=0.9\linewidth]{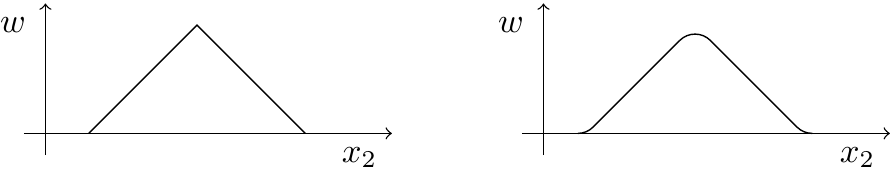}
\caption{Sketch of the effect of the  mollification in ({3})
in a direction orthogonal to the fold.
Left panel: the distance from the boundary $\mathrm{dist}(x,\partial\Omega)$ is a function with slope $\pm 1$
and sharp kinks. Right panel: the mollification defined in ({3})
still has slope $\pm1$ on large parts of the domain, but has smooth transitions from one
value to the other
over a length of the order of $\sigma$, see Figure {3}.}
\end{figure}

\begin{figure}[t]\centering
\includegraphics[width=0.5\linewidth]{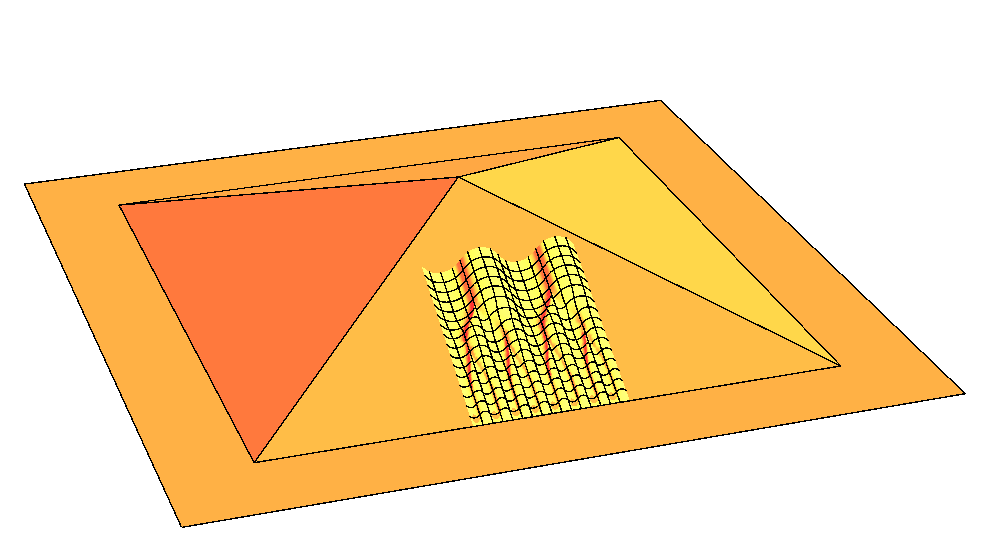}\hfill
\includegraphics[width=0.4\linewidth]{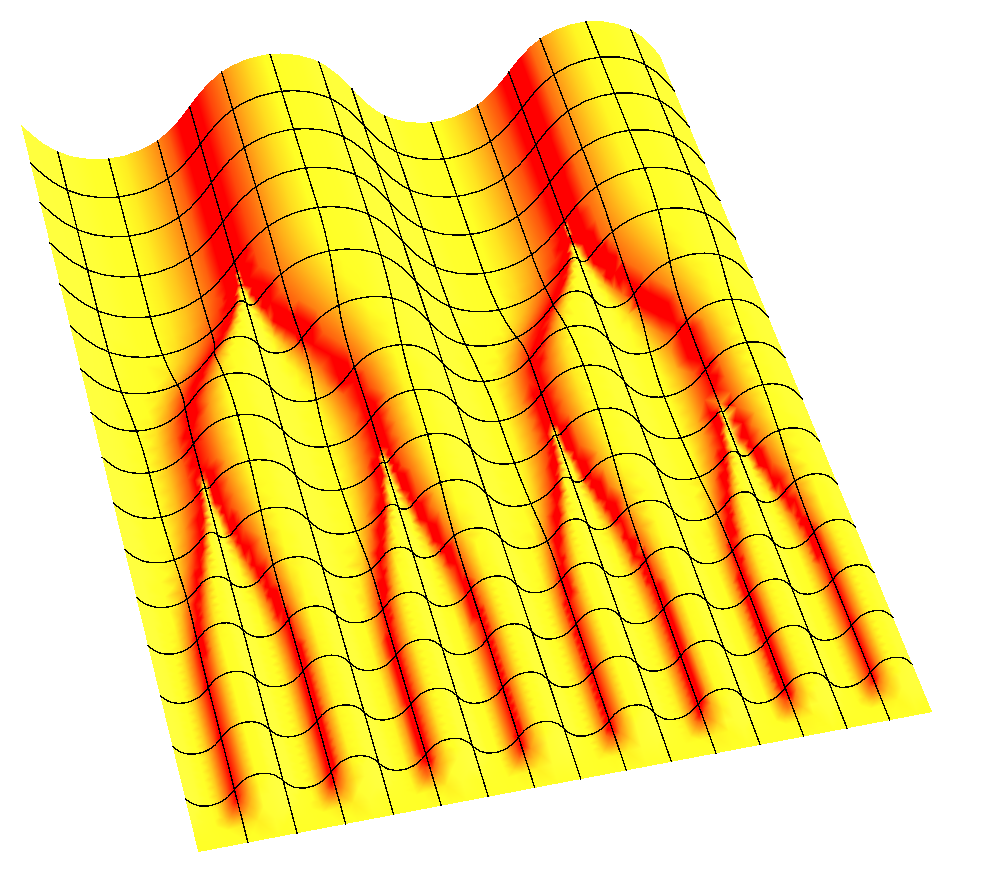}
\caption{Sketch of a deformation achieving the optimal upper bound in ({5}).
As in Figure {1}, the debonded region $\Omega=(0,1)^2$ is a square.
The starting point, at a coarse scale, is the ``tent''-form
illustrated in Figure  {1}. At a finer scale, folds orthogonal to the boundary relax the
tangential compression (left panel, folds are only drawn in a small region). The period of the folds is of order $h$ close to the boundary, and
via a sequence of period-doubling steps becomes coarser in the inside (right, blow-up of the folds from the middle panel).
}
\end{figure}

\section{Pattern formation in debonded thin films}

A finer analysis of the nonlinear elasticity model that had led to ({1}) showed that, in the case of
isotropic compression, also the in-plane components exhibit fine-scale oscillations which refine close to the boundary
{[22, 23, 24, 25]}.
This analysis was based on the F\"oppl-von K\'arm\'an model, which includes the tangential
components of the displacement $u$ as well. After rescaling the energy takes the form (in the case of zero Poisson's ratio for simplicity)
\begin{equation}
  E_{\sigma}[u,w]=\int_\Omega \left( |Du+Du^T+Dw\otimes Dw-\mathrm{Id}|^2
  +\sigma^2 |D^2w|^2 \right)dx\,.
\end{equation}
Here $\Omega\subset\mathbb{R}^2$ is, as above, the debonded region, and the displacements
$u$ and $w$ vanish at the boundary of $\Omega$, corresponding to the fact that the rest of the film
is still bound to the substrate. The isotropic compressive strain has been scaled to 1, and one can check that
$E_\sigma[0,w]=1+I_\sigma[w]$.
The key result from
{[22, 23]} was that the minimum energy scales proportional to $\sigma$:
\begin{theorem}[From {[22, 23]}]
Let $\Omega\subset\mathbb{R}^2$ be a bounded domain with piecewise smooth boundary. Then there are two
constants $c_L,c_U>0$ such that
\begin{equation}
 c_L\sigma \le \min \{E_{\sigma}[u,w]: u=0, w=0 \text{ on $\partial\Omega$} \} \le c_U\sigma\,.
\end{equation}
\end{theorem}
The argument used for proving the lower bound also proves that a
finite fraction of the energy is localized in a thin strip close to the boundary.

Similar statements hold if the plate theory in ({4}) is replaced by a fully
three-dimensional nonlinear elastic model. For $v:\Omega\times(0,h)\to\mathbb{R}^3$, $h>0$, we define
\begin{equation}
 E_{h}^{3D}[v] = \frac1h\int_{\Omega\times(0,h)} W(Dv) dx
\end{equation}
where $W:\mathbb{R}^{3\times3}\to[0,\infty)$ is the elastic stored energy density, which vanishes on
the set of proper rotations $\mathrm{SO}(3)$ and has quadratic growth, in the sense that
\begin{equation}
 c\, \mathrm{dist}^2(F, \mathrm{SO}(3)) \le W(F) \le c' \,\mathrm{dist}^2(F,\mathrm{SO}(3))
\end{equation}
for some positive constants $c$ and $c'$. The factor $1/h$ is included explicitly in ({6})
to obtain an energy per unit thickness, corresponding to ({4}).

In the nonlinear case the thickness $h$ of the film and the
compression $\delta$ enter the problem separately, however to leading order and after scaling the optimal energy only depends on
the combination $\sigma=h/\delta^{1/2}$. In order to understand this expression it is instructive to recall
the relation between the three-dimensional problem  $E_{h}^{3D}$ and its
two-dimensional counterpart $E_{\sigma}$. In particular, a given pair $(u,w)$ in ({4}) corresponds to
a three-dimensional deformation $v_\delta$ of the form
\begin{equation}
v_\delta(x_1,x_2,x_3)=(1-\delta) \left[ \psi(x_1,x_2) + x_3 n(x_1,x_2)\right]
 \end{equation}
 where
  \begin{equation}
\psi(x_1,x_2)=
\begin{pmatrix}
x_1+2\delta u_1(x_1,x_2)\\
x_2+2\delta u_2(x_1,x_2)\\
(2\delta)^{1/2} w(x_1,x_2)
 \end{pmatrix}
 \end{equation}
 represents the deformation of the $x_3=0$ layer
 and
  \begin{equation}
n(x_1,x_2)=
\begin{pmatrix}
-(2\delta)^{1/2}\partial_1 w(x_1,x_2)\\
-(2\delta)^{1/2}\partial_2 w(x_1,x_2)\\
 1
\end{pmatrix}
 \end{equation}
 is, to leading order, the normal to the surface described by $\psi$ and gives the out-of-plane component of the strain.
An expansion of $E_{h}^{3D}[v_\delta]$ for small $\delta$ shows that the leading order
contribution is proportional to $\delta^2 E_{\sigma}[u,w]$  if the Poisson's ratio of the material
vanishes. See for example [22, App. A and App. B] for a more detailed discussion of this point. A rigorous relation between  $ E_{\sigma}$ and
 $E_{h}^{3D}$ was derived in {[26, 27]}
 by means of $\Gamma$-convergence, these results however
 are appropriate for a different regime, with much smaller energy, and therefore
 do not apply directly to the situation of interest here.
\begin{theorem}[From {[25]}]
Let $\Omega\subset\mathbb{R}^2$ be a bounded domain with piecewise smooth boundary, $\delta\in (0,1)$, $h\in (0,\delta^{1/2})$. Then there are two
constants $c_L,c_U>0$ such that
\begin{equation}
 c_L\sigma \le \min \{\frac{1}{\delta^2}E_{h}^{3D}[u]: u(x)=(1-\delta)x \text{ for $(x_1,x_2)\in \partial\Omega$} \} \le c_U\sigma
\end{equation}
where $\sigma=h/\delta^{1/2}$.
\end{theorem}

The significance of Theorem {1} and Theorem {2} is best understood
by considering the key ideas in the proofs.
The upper
bound in ({5}) and ({11})
is proven by explicitly constructing a suitable deformation
field $(u,w)$. This is done in several steps.
The first step is the Ortiz-Gioia construction given in ({3}), which correctly
describes the large-scale behavior of the film and relaxes the compression in direction normal
to the boundary, as in Figure {1}. In the second step one adds fine-scale
oscillations in the orthogonal direction, as illustrated in Figure {3}. This microstructure
does not change the average shape significantly but  relaxes the strain component
tangential to the boundary. Finally, one realizes that optimal deformations have
oscillations on a very fine scale close to the boundary, to adequately match the boundary data,
but much coarser oscillations in the interior, to minimize the bending energy.
Therefore a number of period-doubling steps are inserted, as illustrated in Figure {3}.
Analogous self-similar branched patterns had previously appeared in the study of microstructures
in shape-memory alloys {[28, 29]}, where for a simplified model it had been possible to
show that minimizers are indeed asymptotically self-similar {[30]}.
The scaling in the presence of finite elasticity, both of the martensite and in the surrounding austenite,
was then studied in {[31, 32]}; vectorial variants of the model
were considered in  {[33, 34]}.
A similar approach has been useful also for a variety of other problems, ranging from
magnetic patterns in ferromagnets  {[35, 36, 37]} to field penetration in superconductors
{[38, 39]}, dislocation structures in crystal plasticity {[40]}
and coarsening in thin film growth
{[41]}.

This variational approach to microstructure formation in thin elastic sheets is much more general, and indeed it can be applied to a number of
related problems. One example is paper crumpling {[42, 43]} in which a thin plate, completely detached from
the substrate, is confined to a small volume. In this case it has been possible to construct deformations
with much smaller energy per unit volume. In particular one can obtain an  energy  per unit thickness
proportional to $h^{5/3}$ {[44, 45]}, and one can  approximate
any compressive deformation with this energy.

\begin{theorem}[From {[45]}]
Let $\Omega\subset\mathbb{R}^2$ be a bounded domain, $r>0$. Then there is a map $v:\Omega\times(0,h)\to B_r(0)$ such that
\begin{equation}
 E_{h}^{3D}[v]\le c h^{5/3}\,.
\end{equation}
The constant $c$ may depend on $\Omega$ and $r$ but not on $h$. Further, if $v_0:\Omega\to \mathbb{R}^3$ is a short map, i.e.,
a map which obeys $|v_0(x)-v_0(y)|\le |x-y|$ for all $x,y\in\Omega$, then there is a sequence $v_h$, converging to $v_0$,
such that
\begin{equation}
 \lim_{h\to0}
 \frac1{h^{\alpha}} E_h^{3D}[v_h]=0
\end{equation}
for any $\alpha<5/3$.
Convergence of $v_h$ is understood as uniform convergence of the vertical averages.
\end{theorem}
The proof of this is based on the combination of three ingredients. The first one is an approximation
of short maps with origami maps:
\begin{theorem}[From {[45]}]
Let $v_0:\Omega\to \mathbb{R}^3$ be a short map, i.e.,
a map which obeys $|v_0(x)-v_0(y)|\le |x-y|$ for all $x,y\in\Omega$. Then there is a sequence $v_j$
of Origami maps converging uniformly to $v_0$.
\end{theorem}
Here we say that a map $v:\mathbb{R}^2\to\mathbb{R}^3$ is an Origami map if it is continuous and piecewise isometric,
i.e., if the domain can be subdivided into pieces such that $v$ is a linear isometry (a translation plus a rotation)
in each piece. The number of pieces is allowed to diverge only at infinity, in the sense that  only finitely many pieces are allowed in any bounded
subset of $\mathbb{R}^2$.

The second step is to approximate any  Origami maps with low-energy maps:
\begin{theorem}[From {[45]}]
Let $\Omega\subset\mathbb{R}^2$ be a bounded domain, $v_0:\Omega\times(0,h)\to B_r(0)$ be an Origami map.
Then for any Origami map $v_0$ there is a sequence of maps $v_h:\Omega\times(0,h)\to\mathbb{R}^3$, converging to $v_0$,
such that
\begin{equation}
E_h^{3D}[v_h]\le C h^{5/3}\,.
\end{equation}
The constant may depend on $\Omega$ and $v_0$ but not on $h$.
\end{theorem}
This is proven by an explicit construction around each fold.

Another  related problem of high current interest is the study of wrinkling patterns in  graphene
sheets {[46, 47]}. This has been addressed by a similar model,
in which the boundary conditions are replaced by a viscous term describing the interaction
with a substrate {[48, 49, 50]}. It would be interesting to
see if the methods discussed here can be useful also for this variant of the problem.

\section{Pattern formation in rebonded thin films}
The microstructures spontaneously developed by compressed thin films
can be controlled if the
 geometry of the debonded region is designed appropriately {[5, 6]}.
 One possibility is to introduce a sacrificial layer between
the film and the substrate, and then to selectively etch away a part of it, so that the boundary of the debonded region is straight,
see sketch in Figure {4}.
The film then partially rebonds to the surface, leading to complex patterns of tunnels,
which in some cases refine close to the boundary, see Figure {5}.

These patterns can be studied by coupling the von-K\'arm\'an energy with a fracture term
 proportional to the debonded area,
\begin{equation}
  E_{\sigma,\gamma}[u,w]=\int_\Omega \left( |Du+Du^T+Dw\otimes Dw-\mathrm{Id}|^2 +\sigma^2 |D^2w|^2 \right)dx
  + \gamma |\{w>0\}|.
\end{equation}
The three terms represent stretching, bending and bonding energies respectively.
Here $u:\Omega\to\mathbb{R}^2$ are the (scaled) tangential displacements and $\gamma>0$ is the
bonding energy per unit area (related to Griffith's fracture energy),
$|\{w>0\}|$ represents the area of the set where the vertical displacement $w$ is nonzero. Equivalently one could take the debonded
state as reference and consider a negative term proportional to the rebonded area,
$-\gamma|\{w=0\}|$; the two energies only differ by an additive constant.
The appropriate boundary conditions correspond to the film being bound to a substrate on one side
of the domain; for simplicity we shall focus on $\Omega=(0,1)^2$ with $u=0$ and $w=0$ on the
$\{x_1=0\}$ side of $\Omega$. As above, we assume $w\ge 0$ everywhere.

\begin{figure}[t]\centering
\includegraphics[width=0.6\linewidth]{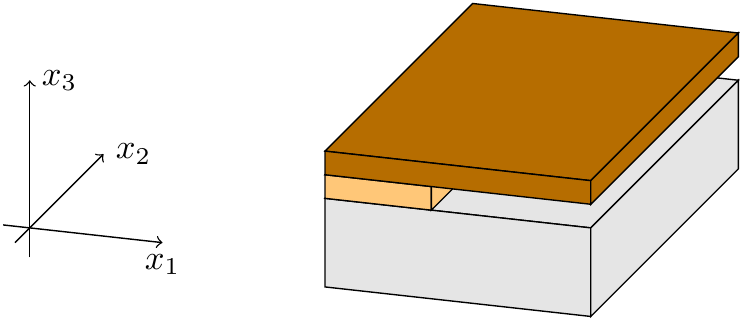}\hfill
\caption{Geometry of the partially delaminated film.
The intermediate sacrificial layer is removed chemically only for $x_1>0$.
The free-standing film is subject to compression at the Dirichlet boundary
$x_1=0$ and may rebond to the substrate. }
\end{figure}

The mathematical analysis of the  energy
({15}) leads to the rich phase diagram sketched in Figure {6}, which contains four different regimes {[51]}
that we now illustrate.

For large specific bonding energy $\gamma$ the film is completely bound to the substrate.
In particular the film is flat, so that there is no
bending energy, but the stretching energy is not released. The total energy is then proportional to the area of $\Omega$, and one obtains
$E_{\sigma,\gamma}[0,0]=2$. This is regime A in Figure {6} and Theorem {6}.

\begin{figure}[t]\centering
\caption{Experimental picture of tube branching in
 Si$_{1-x}$Ge$_x$ film on a thick SiO$_2$ substrate.
 Experimental picture
removed for
copyright reasons, please see [5, Fig. 2].}
\end{figure}

The opposite case of  very small bonding energy $\gamma$ is also easy to understand after the foregoing discussion: here the bonding term plays no significant role and the
film is completely detached from the substrate. One recovers the result of the blistering problem
of Theorem {1},
$E_{\sigma,\gamma}[u,w]\simeq E_\sigma[u,w]\le c \sigma$. The corresponding deformations are those illustrated in Figure {3}.
This is regime D in Figure {6} and Theorem {6}.

\begin{figure}[t]\centering
\includegraphics[height=0.6\linewidth]{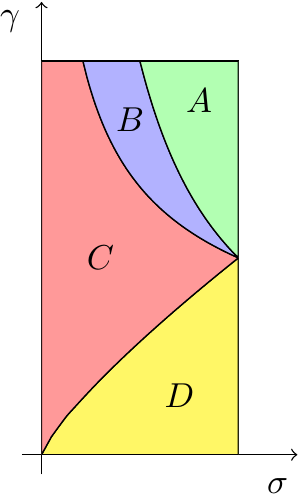}\hfill
\caption{Phase diagram for $E_{\sigma,\gamma}[u,w]$ in the $(\sigma,\gamma)$ plane.}
\end{figure}

For intermediate values of $\gamma$ the situation is more complex, in particular
debonded channels are formed, which separate wider bonded regions.
In regime B the pattern is periodic and, away from the Dirichlet boundary,
depends only on the tangential variable $x_2$.
A large part of the film is bonded to the substrate, but bonded regions are separated
by thin tubes, see Figure {7}.
Denoting by $h$ the period of the oscillations, and by $\delta$ the width of a tube, the
total volume fraction of the tubes is $\delta/h$, therefore the bonding energy is proportional
to $\gamma\delta/h$. Each tube has to release a compression of $h$ over a width $\delta$, therefore
the term $|Dw|^2$ is of order $h/\delta$ inside the tubes (the stretching energy is then completely relaxed). This gives
$|Dw|\sim (h/\delta)^{1/2}$ in the tubes, and hence $|D^2w|\sim (h/\delta)^{1/2}/\delta$. Therefore
the total energy can be estimated by
\begin{equation}
 \gamma \frac{\delta}{h}+
 \sigma^2 \frac{\delta}{h} \left( \frac{h^{1/2}/\delta^{1/2}}{\delta}\right)^2 =
 \gamma \frac{\delta}{h}+
  \frac{\sigma^2}{\delta^2}\,.
\end{equation}
Optimizing in $\delta$ we obtain $\delta\sim  \sigma^{2/3}h^{1/3}\gamma^{-1/3}$
(this is clearly only admissible if $\delta\le h\le 1$).
The
period $h$ is fixed by the energetic cost of the interpolation region close to the boundary. In this part of the domain
there is no stretch-free construction, and indeed an interpolation over a boundary layer of thickness
$\varepsilon$ results in a total stretching energy of $\varepsilon(1+ h^2/\varepsilon^2)$. Optimizing over $\varepsilon$
we obtain $\varepsilon\sim h$, and therefore the total energy for the laminate construction is
\begin{equation}
h +
\gamma \frac{\delta}{h}+
  \frac{\sigma^2}{\delta^2}\,.
\end{equation}
Inserting the value of $\delta$ obtained above and minimizing in $h$ we conclude that
$h$ and $E$ are proportional to $(\sigma\gamma)^{2/5}$.
The width of each tube $\delta$ is then proportional to
$\sigma^{4/5}\gamma^{-1/5}$. This is regime B in Figure {6};
a precise version of this construction proves the second bound in Theorem {6}.

If the bending term becomes more important, it is convenient to insert period-doubling steps,
just like in the discussion of the functional ({4}).
The resulting pattern is shown in   Figure {8}. In comparison to the pattern
of Figure {3} the key difference is that the bending is localized to a small region, whereas
large parts of the film are bond to the substrate. The period-doubling steps are only possible at the
expense of stretching energy; balancing the different terms one finds {[51]} that the resulting energy is proportional to
$\sigma^{1/2}\gamma^{5/8}$.
The result of the construction is summarized in the following statement.

\begin{figure}[t]\centering
\includegraphics[width=0.95\linewidth]{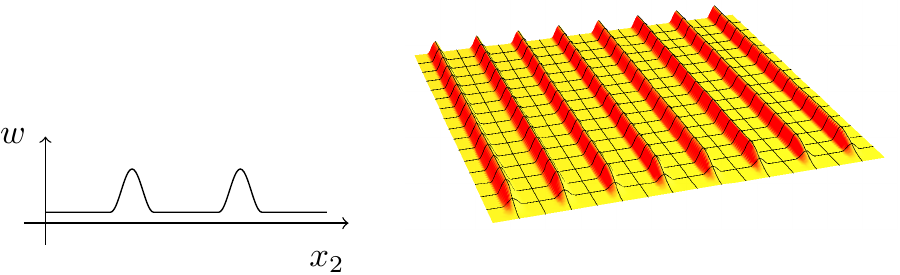}\hfill
\caption{Sketch of the laminate regime (B).}
\end{figure}

\begin{figure}[t]\centering
\includegraphics[width=0.75\linewidth]{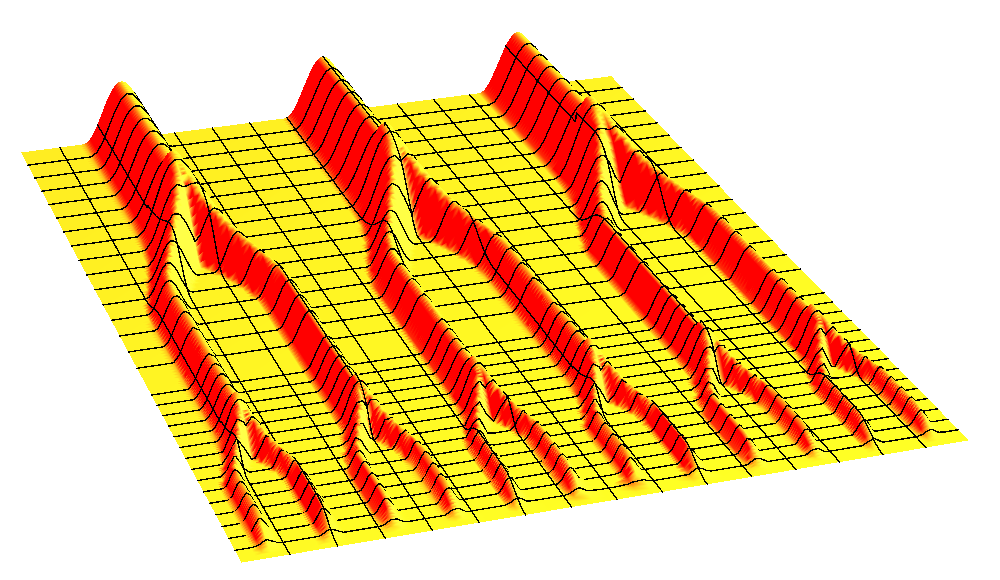}
\caption{Sketch of the tube branching regime (C).}
\end{figure}

\begin{theorem}[From {[51]}]
 Let $\gamma>0$, $\sigma\in(0,1)$. There are  $u$, $w$ which obey the stated boundary conditions and
  \begin{equation}
   E_{\sigma,\gamma}[u,w]\le c\left\{\!\!\!
  \begin{tabular}{lll}
    $1$ &  if  $\sigma\gamma>1$ & \text{(regime A),}\\
    $(\sigma\gamma)^{2/5}$ &  if  $\sigma^{-4/9}\le \gamma\le \sigma^{-1}$
    & \text{(regime B),}\\
    $\sigma^{1/2}\gamma^{5/8}$ & if  $\sigma^{4/5}\le \gamma\le
    \sigma^{-4/9}$ & \text{(regime C),}\\
    $\sigma$ &  if $\gamma<\sigma^{4/5}$ & \text{(regime D).}
  \end{tabular}\right.
\end{equation}
\end{theorem}
The proof is based on making the constructions sketched above precise,
details are given in {[22]} for regime D and in {[51]}
for regimes B and C. Regime A, as discussed above, is immediate.


Optimality of the phase diagram just discussed can be at least partially proven by providing matching lower bounds on the energy.
In particular, one can show the following.
\begin{theorem}[From {[51]}]
 Let $\gamma>0$, $\sigma\in(0,1)$. For any $u$, $w$ which obey the stated boundary conditions one has
 \begin{equation}
  E_{\gamma,\sigma}[u,w]\ge c
\left\{\!\!\!
  \begin{tabular}{lll}
    $1$ &  if  $\sigma\gamma>1$ & \text{(regime A),}\\
    $(\sigma\gamma)^{2/3}$ &  if  $\sigma^{1/2}\le \gamma\le \sigma^{-1}$
     & \text{(regime B'),}\\
    $\sigma$ &  if $\gamma<\sigma^{1/2}$ & \text{(regime D').}
  \end{tabular}\right.
 \end{equation}
\end{theorem}
Whereas the statement in regime D' follows from {[22]}, the other two bounds are proven in
{[51]} using the Korn-Poincar\'e inequality for $SBD^2$ functions obtained in
{[52]}.

Theorem {7} proves optimality in phases A and D. The bound in the intermediate region does not, however, match the upper bounds stated in Theorem
{6}. Therefore it is at this stage not clear if the branching patterns illustrated in Figure {8} are optimal.

\section{Linear stability analysis}
The general form of the linearized F\"oppl-von K\'arm\'an plate theory under isotropic
compression is {[53, 54]}
\begin{equation}
  E_{\mathrm{FvK}}[u, w] = \frac12 Y h \int_\Omega \left[
(1-\nu) |\epsilon|^2 + \nu (\mathrm{Tr}\epsilon)^2
  + \frac{h^2}{12}
\left[ (1-\nu)  |D^2 w|^2 + \nu (\Delta w)^2\right] \right]dx\,,
\end{equation}
see also {[2, 22]} for a discussion in the present context and
{[27]} for a rigorous mathematical derivation.
Here $\nu\in[-1,1/2]$ is the Poisson ratio, $Y$ Young's modulus, $h$ the film thickness, and
the strain  $\epsilon$ is defined by
\begin{equation}
  \epsilon =  D u + (D u)^T
  + Dw \otimes D w - 2\delta \mathrm{Id}\,,
\end{equation}
where $\delta$ is the eigenstrain (i.e., the compression enforced by the substrate).
We recall that  we use $|M|^2=\mathrm{Tr} M^TM$ for the matrix norm.
For $\nu=0$, after a rescaling ({20}) reduces to ({4}).
We recall that in  [22, App. B] it was shown that the scaling behavior of the functional 
$ E_{\mathrm{FvK}}$ is the same for all $\nu \in (-1,1/2]$, hence our results hold also for generic
values of the Poisson ratio. Of course, the regime $\nu\ge0$ is the most relevant.

For small $\delta$ one can linearize around the state $u=0$, $w=0$. After straightforward
computations this leads to
\begin{alignat*}1
  E_{\mathrm{FvK}}^\mathrm{lin}[u, w] = \frac12 Y h \int_\Omega &
\Big[(1-\nu) |Du+Du^T -2\delta\mathrm{Id}|^2
-4\delta |Dw|^2
+ \nu (2\mathrm{div } u-4\delta)^2\\
&  -8\delta \nu  |Dw|^2
  + \frac{h^2}{12}
\left[ (1-\nu)  |D^2 w|^2 + \nu (\Delta w)^2 \right]\Big]dx\,.
\end{alignat*}
In this linearized functional $u$ and $w$ are decoupled. The dependence on $u$ in convex, hence
$u=0$ is the minimizer with the given boundary data. The dependence on $w$ is however not necessarily
convex. Working for concreteness in a circle of radius $R$, we can assume $w$ to be radial, $w(x)=\varphi(|x|)$,
subject to $\varphi(R)=0$, so that
\begin{equation*}
 Dw(x)=\varphi'(|x|) \frac{x}{|x|}
\end{equation*}
and
\begin{equation*}
 D^2w(x)=\varphi''(|x|) \frac{x}{|x|}\otimes \frac{x}{|x|}
 + \varphi'(|x|) \left(\frac{\mathrm{Id}}{|x|}-\frac{x\otimes x}{|x|^3}\right)\,.
\end{equation*}
Inserting into the energy leads to the one-dimensional variational problem
\begin{alignat*}1
\frac12 Y h\int_0^R &\left[
-4\delta (1+2\nu) (\varphi'(r))^2
\phantom{\left(\frac{\varphi'(r)}r\right)^2}\right.\\
&\left.\hskip1cm
+ \frac{h^2}{12}
\left[  (\varphi''(r))^2 +
\left(\frac{\varphi'(r)}r\right)^2
+2\nu \frac{\varphi'(r)\varphi''(r)}r\right] \right] r dr\,.
\end{alignat*}
This is positive definite if the first term, of order $\delta$, is not larger then the second
term, of order $h^2/R^2$. Therefore the loss of stability, which corresponds to Euler buckling, occurs at strains
$\delta\sim h^2/R^2$.
Inserting the experimental data from {[5]}, namely,
$h\sim 20$ nm, $R\sim 10\, \mu$m, $\nu \sim 0.277$, leads to
$\delta_{\mathrm{crit}}\sim 4\cdot 10^{-6}$, which corresponds to a strain of $0.0004$\%.
This is over three orders of magnitude smaller than the experimentally applied strain
$\delta_{\mathrm{Exp}}\sim 0.011=1.1\%$.
Therefore the experiments we discussed take place well beyond the loss of linear stability, and a
buckling-postbuckling analysis does not seem appropriate to understand the deformations. Our variational
approach is instead constructed to deal with deformations and microstructures that appear in the deeply nonlinear regime
and is therefore more suitable to study the mentioned experiments.

\section*{Acknowledgements}
This work was partially supported by the Deutsche Forschungsgemeinschaft
through SFB 1060.



\end{document}